\documentstyle[12pt]{article}
\begin{document}
\newtheorem{thm}{Theorem}[section]
\newtheorem{lem}{Lemma}[section]
\newtheorem{rem}{Remark}[section]
\newtheorem{prop}{Proposition}[section]
\newtheorem{cor}{Corollary}[section]
\newtheorem{defn}{Definition}[section]
\title
{Some inequalities related to isoperimetric inequalities with partial free boundary\\}
\author  {
Meijun Zhu\\
Department of Mathematics\\
 University of Oklahoma\\
Norman, Oklahoma 73019 \\
E-mail: mzhu@math.ou.edu
}

\newcommand{\A} {\alpha}
\date{ }
\maketitle
\input { amssym.def}
\input { amssym.def}
\pagenumbering{arabic}
\newcommand{\bg} {\begin{equation}}
\newcommand{\ede} {\end{equation}}
\newcommand{\M} {\Omega}
\newcommand{\E}{\lambda}
\newtheorem{definition}{Definition}[section]

\maketitle
\makeatletter
\@addtoreset{equation}{section}
\def\theequation{\thesection.\arabic{equation}}
\makeatother
\makeatletter

\maketitle
\begin{abstract}

The main purpose of this paper is to prove a sharp Sobolev inequality  in the exterior of a convex bounded domain. There are two ingredients in the proof:  One is   the observation of   some  new isoperimetric inequalities with partial free boundary,  and the other is an integral inequality (due to Duff \cite{Du}) for any nonnegative function under  Schwarz equimeasurable rearrangement. These ingredients also allow us to establish some Moser-Trudinger type inequalities, and obtain some estimates on the principal frequency of a membrane with partial free boundary, which extend early results of Nehari \cite{Ne} and Bandle \cite{Ba} for two dimensional domains to the one for  any dimensional domains (dimension $ \ge 2$).
\end{abstract}

\section{Introduction}

One of the main motivations of our study is to prove the following sharp Sobolev inequality.

\begin{thm} Suppose that $\M$ is a bounded smooth convex  domain in $\Bbb{R}^n$ ($n \ge 2$) and $\M^c=\Bbb{R}^n \setminus \overline \M$ is the exterior of $\M$. For $1 \le p <n$, we denote $p_*=np/(n-p)$, and $k(n,p)=\pi^{-1/2} n^{-1/p}(\frac {p-1}{n-p})^{1-1/p} \{ \frac {\Gamma(1+n/2) \Gamma(n)}{\Gamma(n/p) \Gamma(1+n-n/p)} \}^{1/n}$ for $p>1$
and  $k(n,1) = \lim_{p \to 1+} k(n, p)=\pi^{-1/2} n^{-1} \{\Gamma(1+n/2)\}^{1/n}$ as the best constants of Sobolev inequalities on $\Bbb{R}^n$. Then
\bg
 (\int_{\M^c} |u|^{p_*} )^{1/p_*}  \le  2^{1/n}k(n,p)  (\int_{{ \M^c}}
|\nabla u|^p)^{1/p} 
\label{A1}
\ede
holds for all $ u\in C_0^\infty (\Bbb{R}^n)$.
\label{thmA1}
\end{thm}

In the rest of the paper, we shall keep the same notations for $n, \ p, \ p_*$  and $  k(n,p).$

\medskip

The study of  the sharp Sobolev inequalities goes back to  Aubin and Talenti. In 1976, Aubin \cite{au}, Talenti \cite{Ta} respectively, obtained  the sharp Sobolev inequality on $\Bbb{R}^n$ ($n \ge 2$):  Inequality
\begin{equation}
 (\int_{\Bbb{R}^n} |u|^{p_*} )^{\frac 1{p_*}}  \le k(n,p)  (\int_{\Bbb{R}^n}
|\nabla u|^p)^{1/p}
\label{A2}
\end{equation}
holds  for all $u(x)$ satisfying $u\in L^{p_*}(\Bbb{R}^n)$  and $\nabla u\in  L^{p}(\Bbb{R}^n)$.
For $1<p<n$, extreme functions of (\ref{A2})  were also found by them.

Using an even reflection, we easily obtain the sharp  Sobolev inequalities on the upper half space $\Bbb{R}^n_+$: If $u(x)$ satisfies $u\in L^{p_*}(\Bbb{R}^n_+)$  and $\nabla u\in  L^{p}(\Bbb{R}^n_+)$, then 
$$
 (\int_{\Bbb{R}^n_+} |u|^{p_*} )^{\frac 1{p_*}}  \le 2^{1/n}k(n,p)  (\int_{\Bbb{R}^n}
|\nabla u|^p)^{1/p}.
$$

For $p=2$ (and temporarily assuming that $n\ge 3$), noticing that the upper half space is conformally equivalent to a ball $B_r(0)$  with radius $r$, we can write, using the conformal invariance of the conformal Laplacian operator, the sharp Sobolev inequalities on $B_r(0)$ and on the exterior of $B_r(0)$ as follows:
\begin{equation}
 (\int_{B_r(0)} |u|^{\frac{2n}{n-2}} )^{\frac{n-2}n}  \le 2^{2/n} k^2(n,2)  \bigg(\int_{B_r(0)} |\nabla u|^2+ \frac {n-2}{2r}\int_{\partial B_r(0)}|u|^2 \bigg), \ \ \ \ \ \ \ \ \forall u \in H^1(\Omega)
\label{A4}
\end{equation}
and
\begin{equation}
 (\int_{B^c_r(0)} |u|^{\frac{2n}{n-2}} )^{\frac{n-2}n}  \le 2^{2/n} k^2(n,2)  \bigg(\int_{B^c_r(0)} |\nabla u|^2- \frac {n-2}{2r}\int_{\partial B_r(0)}|u|^2 \bigg), \ \ \ \ \ \ \ \ \forall u \in C_0^\infty(\Bbb{R}^n).
\label{A5}
\end{equation}
 Naturally (parallel to Aubin's conjecture which were proved by Hebey and Vaugon in \cite{HV}), one may ask if   there is  any sharp Sobolev inequality for any smooth bounded domain $\Omega$.  Namely, is there any constant $C(\Omega)$ so that
\begin{equation}
 (\int_\Omega |u|^{\frac{2n}{n-2}} )^{\frac{n-2}n}  \le  2^{2/n} k^2(n,2) \int_{\Omega}
|\nabla u|^2  + C(\Omega) \int_{\partial \Omega}u^2
\label{A6}
\end{equation}
holds for all $u \in H^1(\Omega)$?

Such an existence result of $C(\M)$ was proved by the author in a  joint work with Y.Y. Li \cite{LZgafa}. Recently there are many results concerning some  similar inequalities, see for example,  Hebey and Vaugon \cite{HV},  Hebey \cite{H},  Y.Y.Li and M.Zhu \cite{LZgafa}-\cite{LZ}, Zhu \cite{Zhu1}-\cite{Zhu2}, Druet \cite{Dr}, Aubin, Druet and Hebey \cite{ADH}, Aubin and Li \cite{AL},  and references therein. The uniform way to prove these inequalities is to argue by contradiction, which yields no information on $C(\M)$. So far, to my knowledge, there is no any upper bound estimate on $C(\M)$ unless  $\partial \M$ is close to a sphere (see, for example, Pan and Wang \cite{PW}). From analytic point of view, it remains as  a challenge to determine the optimal constant for $C(\M)$.  Estimates on the optimal constant  will  shed light on the solvability of some elliptic equations involving critical Sobolev exponents, and amazingly will have some geometric impacts (especially on isoperimetric inequalities) as we will see below.

 Let us denote $C_{opt}(\Omega)$ as the optimal constant for $C(\M)$ in (\ref{A6}).  It is well known that $C_{opt}(\Omega)$ is closely  related to the geometric property of the boundary $\partial \Omega$. For instance, due to the work of Wang \cite{W} (see also, Adimurthi and Mancini \cite{AM}), we know that for any bounded smooth domain $\Omega$, $C_{opt}(\Omega) \ge (n-2)2^{2/n-1}  k^2(n,2) \cdot \max_{x \in \partial \Omega} H(x)$, where $H(x)$ is the mean curvature function of $\partial \Omega$ with respect to the inner normal of $\partial \Omega$ (e.g., the unit sphere in $\Bbb{R}^n$ has positive mean curvature).  Nevertheless, in view of (\ref{A5}), we may ask: Given a bounded convex smooth domain $\M$ in $\Bbb{R}^n$ (assuming $n \ge 3$), does the inequality
$$
(\int_{ \Omega^c} |u|^{\frac{2n}{n-2}} )^{\frac{n-2}n}  \le 2^{2/n}  k^2(n,2)  \int_{{ \Omega^c}}
|\nabla u|^2
$$
hold for all $u \in C_0^\infty(\Bbb{R}^n)$?

Obviously, this is a special case of the following question.

\noindent{\bf  Question 1.1.}\ Given a bounded smooth convex domain $\M$ in $\Bbb{R}^n$ (assuming $n \ge 2$), Does the inequality
$$
 (\int_{ \Omega^c} |u|^{p_*} )^{1/p_*}  \le  2^{1/n}k(n,p)  (\int_{{ \Omega^c}}
|\nabla u|^p)^{1/p} 
$$
hold for all $ u \in C_0^\infty(\Bbb{R}^n)$?

\medskip

Our Theorem \ref{thmA1} gives an affirmative answer to this question.

\medskip

One of the most interesting consequences (however, it is a {\it fake} consequence) of Theorem \ref{thmA1} would be the following ``isoperimetric inequality with partial free boundary''. 

\begin{thm}({\bf isoperimetric inequality with partial free boundary})\  Let $\M$ be a bounded piecewise smooth domain in $\Bbb{R}^n$ ($n \ge 2$). Suppose that the boundary $\partial \M$ consists of two  smooth  hypersurfaces $\Gamma_1$ and $\Gamma_2$. If $\Gamma_2$ is concave with respect to $\M$, then
\bg
 \frac{|\Gamma_1|}{|\Omega|^{1-1/n}} \ge \frac { {\pi}^{1/2} n}{\{2 \Gamma(1+n/2)\}^{1/n}},
\label{A7}
\ede
where $|\Gamma_1|$ is the $n-1$-dimensional surface area of $\Gamma_1$ and $|\Omega|$ is the $n$-dimensional volume of $\Omega$.
\label{thmA2}
\end{thm}

Here, we shall specify the concavity of $\Gamma_2$: We say that  $\Gamma_2$ is concave with respect to $\M$ if any line segment between two endpoints on $\Gamma_2$ does not belong to $\Omega$.

Theorem \ref{thmA2} is a consequence of Theorem \ref{thmA1} since it is well known that Theorem \ref{thmA1} implies Theorem \ref{thmA2} (by choosing $p=1$ and $u(x)$ being the approximation of a characteristic function). However, we say that this is a {\it fake} consequence since,  in this paper, we will use   Theorem \ref{thmA2} to prove  Theorem \ref{thmA1}  via Schwarz symmetrization.  Actually this is the main idea in the proof of Theorem \ref{thmA1}. It turns out that the proof of Theorem \ref{thmA2} becomes crucial.

We prove Theorem \ref{thmA2} by using some ideas similar to Steiner symmetrization. To distinguish our method, we shall call it ``the method of reflection Steiner symmetrization''. The details  will be addressed in Section 2. 

In Section 3  we discuss the procedure of equimeasurable rearrangement (Schwarz symmetrization) and reprove an integral inequality which was initially given by Duff \cite{Du}. Such an  integral inequality allows us to prove Theorem \ref{thmA1} by using Theorem \ref{thmA2} in Section 4. Actually, in Section 4  we prove the following  more general sharp Sobolev inequality:

\begin{thm}
{\bf (sharp Sobolev inequality)}\ Let $\M$ be a bounded piecewise smooth domain in $\Bbb{R}^n$ ($n \ge 2$). Suppose that the boundary $\partial \M$ consists of one set of disjoint smooth hypersurfaces  $\Gamma_1$   and another  smooth hypersurface $\Gamma_2$. If $\Gamma_2$ is concave with respect to $\M$, then 
\bg
 (\int_{ \Omega} |u|^{p_*} )^{1/p_*}  \le  2^{1/n}k(n,p)  (\int_{{ \Omega}}
|\nabla u|^p)^{1/p} 
\label{A8}
\ede
holds for all $ u $ satisfying $u\in W^{1,p}(\Omega)$ and $u=0$ on every hypersurface in $\Gamma_1$.
\label{thmA3}
\end{thm}

The crucial lemma in the proof of Theorem \ref{thmA3} is Lemma \ref{lem4-1}. This lemma and other ingredients in the proof of Theorem \ref{thmA3} enable us to obtain other two results.

\medskip

In the remainder of Section 4, first of all, we  establish some Moser-Trudinger type inequalities. Let us start with  a Moser-Trudinger type inequality on a bounded smooth domain $\M$ in $\Bbb {R}^2$.  The following  inequality was proved by Cherrier \cite{Che} (see also Chang and Yang \cite{CP},  more general result concerning piecewise smooth domains was established in \cite{CP}): For any $u \in C^1(\M)$, if $\int_\M |\nabla u|^2 dz \le 1$ and $\int_\M u dz=0$, then
$$
\int_\M e^{ 2 \pi u^2} dz \le C(|\M|),
$$
where $C(|\M|)$ is a constant depending only on the measure of $\M$. This is similar to a ``Moser-Trudinger inequality with boundary'', which we  shall describe below. The standard  Moser-Trudinger inequality says: If $u \in C_0^1(\M)$ and $\int_\M |\nabla u|^2 dz \le 1$, then
$$
\int_\M e^{ 4 \pi u^2} dz \le C_1(|\M|),
$$
where $C_1(|\M|)$ is another constant depending only on the measure of $\M$. 
If $\M$ is a piecewise smooth bounded domain whose boundary consists of two curves $\Gamma_1$ and $\Gamma_2$, where $\Gamma_2$ is part of a straight line, then for any $u \in C^1(\M)$, $u=0$ on $ \Gamma_1$ and $\int_{\M} |\nabla u|^2 dz \le 1$, using the  even reflection with respect to $\Gamma_2$, we have
$$
\int_{\M} e^{ 2 \pi u^2} dz \le C_2(|\M|)
$$
for a constant $C_2(|\M|)$ depending only on the measure of $\M$.

This leads us  to ask: in the inequality (and higher dimensional analogs) due to Cherrier, Chang and Yang, for which kind of domains  one can replace the condition $\int_\M u dz=0$ by assuming that $u$ vanishes on part of the boundary?

We give a complete answer to this  question in the following theorem.

\begin{thm} Let $\M \subset \Bbb{R}^n$ ($n \ge 2$) be a piecewise smooth bounded domain whose boundary consists of one set of disjoint smooth hypersurfaces $\Gamma_1$ and  one smooth  hypersurfaces $\Gamma_2$. If $\Gamma_2$ is concave with respect to $\Omega$, then for any $u(x)$ satisfying $u\in C^1(\Omega)\cap C(\overline \M)$, $\int_{\Omega}|\nabla u|^n \le 1$ and $u=0$ on every hypersurface in $\Gamma_1$, the following inequality holds:
\bg
\int_\M e^{ \beta_n  u^{\frac n{n-1}}} dz \le C_*(|\M|),
\label{A9}
\ede
where $C_*(|\M|)$ is a constant depending only on the measure of $\M$, $\beta_n=n (\frac {\omega_{n-1}} 2)^{1/(n-1)},$ and $\omega_{n-1} $ is the surface area of the $(n-1)$-dimensional  unit sphere.
\label{thmA4}
\end{thm}

The concavity assumption on  $\Gamma_2$ is sharp. We have an example which shows that Theorem \ref{thmA4} does not hold if the concavity assumption is removed.

\medskip

At the end of Section 4, we give some estimates on the principal frequency of a membrane with partial free boundary. Let $\M \in \Bbb{R}^2$ be a bounded piecewise  smooth domain whose boundary $\partial \M$ containing a curve $\alpha$. If $\alpha$ is concave with respect to $\M$ and $\M$ is simply connected, Nehari \cite{Ne} proved 
$$
\Lambda \ge \Lambda_0,
$$
where $\Lambda$ is the principle frequency of a homogeneous membrane whose boundary is free along $\alpha$ and fixed along $\partial \M \setminus \alpha$, and $\Lambda_0$ is the principle frequency of a homogeneous semi-circular membrane of equal measure whose boundary is free along the diameter  and fixed along the semi-circle. Later, Bandle \cite{Ba} was able to remove  the simply connected  assumption  on $\M$. Their methods fail in the higher dimensional cases. Using our new  isoperimetric inequalities, we are able to extend such an estimate to the one for any dimensional domains.

\begin{thm}  Let $\M$ be a bounded piecewise smooth domain in $\Bbb{R}^n$ ($n \ge 2$). Suppose that the boundary $\partial \M$ consists of one set of disjoint smooth hypersurfaces $\Gamma_1$ and another  smooth hypersurface $\Gamma_2$.   If  $\Gamma_2$ is  concave with respect to $\M$, then 
$$
\inf_{u\in H^1(\M) \setminus \{0\}, u=0 \ \mbox{on}\  \Gamma_1} \frac{\int_{\M} |\nabla u|^2 dx}{\int_\M u^2 dx} \ge \inf_{u\in H^1(B_r^+) \setminus \{0\}, u=0 \ \mbox{on}\ \partial B_r^+\setminus \{y=0\}} \frac{\int_{B_r^+} |\nabla u|^2 dx}{\int_{B_r^+} u^2 dx},
$$
where $B_r^+$ is the upper half ball with the center at $\{0\}$  and the same volume as $\M$.
\label{thmA5}
\end{thm}

\noindent{\bf Acknowledgments:} It is my pleasure to thank H. Brezis, X. Cabre, H.D. Cao, S.Y. Cheng, Z-C. Han,  P. Petersen, S.T.Yau, and some of my colleagues in OU for some interesting discussions and their interests in this work.

\section{Isoperimetric inequalities with partial free boundary}
We prove Theorem \ref{thmA2} in this section.

We first observe that (\ref{A7}) holds for $n=2$.
Let $\M$ be a bounded piecewise smooth domain in $\Bbb{R}^2$ whose boundary consists of two smooth curves $\Gamma_1$ and $\Gamma_2$, and  $\Gamma_2$ is concave with respect to $\M$. We can assume that $\Gamma_2$ is not closed, otherwise  the  standard isoperimetric inequality yields $|\Gamma_1|/|\M|^{1/2} \ge 2 \sqrt {\pi} >\sqrt {2 \pi}$. Further, without loss of generality (due to the concavity assumption on $\Gamma_2$), we can assume that $\Gamma_2$ consists of  straight line segments. If $\Gamma_2$ consists of only one straight line segment and $\M$ is on one side of $\Gamma_2$,  by  the  standard isoperimetric inequality and the even reflection with respect to $\Gamma_2$, we have
$$
\frac {|\Gamma_1|} {|\M|^{1/2}} \ge \frac{|\partial B_1^+ \cap \{(x_1, x_2) \ : \ x_2>0\}|}{|B_1^+|^{1/2}} =\sqrt{2 \pi}.
$$
Generally, we denote $\M= \cup_{i=1}^{k} D_i$ where each $D_i$ is on one side of each extended line segment of $\Gamma_2$, and have
$$
 \frac {|\Gamma_1\cap \overline D_i|} {|D_i|^{1/2}} \ge \sqrt{2 \pi}, \ \ \ \ \mbox{for} \ \ \ i=1,2,...k,
$$
which again yields
$$
 \frac {|\Gamma_1|} {|\M|^{1/2}} \ge \frac {\sum_{i=1}^{k}|\Gamma_1 \cap \overline D_i|} {\sum_{i=1}^k |D_i|^{1/2}} \ge \sqrt{2 \pi}.
$$

\medskip

For $n \ge 3$, if $\Gamma_1 \cap \Gamma_2$  lies in the same hyperplane, we can prove (\ref{A7}) by an even reflection just as above. In the rest of this section, we assume that $\Gamma_1 \cap \Gamma_2$ does not lie in the same hyperplane. We are going to symmetrize the domain via some even reflections.  Even though the method is similar to Steiner symmetrization (see, for example, P\'olya and Szeg\"o \cite{PS}), but there are some differences in the procedures of symmetrizations. We shall call our method ``the method of reflection Steiner symmetrization''. 

We argue  by contradiction.
Suppose that for  some  $\epsilon<\frac {1+2^{-1/n}-2^{1-1/n}}{2}$,
\bg
 \frac{|\Gamma_1|}{|\Omega|^{1-1/n}} <(2^{-1/n}-\epsilon) \cdot \frac { {\pi}^{1/2} n}{\{ \Gamma(1+n/2)\}^{1/n}}.
\label{A2-1}
\ede
We are going to derive a contradiction.

We choose a suitable coordinate system so that the area of the domain on the hyperplane $\{(x_1, ..., x_n) \ : \ x_n=0\}$ which is  bounded by the projection from $\Gamma_1 \cap \Gamma_2$ to hyperplane $\{(x_1, ..., x_n) \ : \ x_n=0\}$ is greater than zero, and  $\Gamma_1 \cap \Gamma_2$ is not above (in positive $x_n$-direction) the tangent plane of $\Gamma_2$ which is parallel to $\{(x_1,...,x_n) \ : \ x_n=0\}$.

In the rest of this section, we will always denote $H_\theta$ as a hyperplane which is  parallel to $x_n$-axis and has $\theta$ angle with $x_1$-axis. We need the following two lemmas.

\begin{lem} If $H_\theta$ divides $\M$ into two equal volume domains $\M_1$ and $\M_2$, then at least on  one of them, say $\M_1$,  the following inequality holds:
$$
\frac{|\Gamma_1^{1,1}|}{|\Omega_1|^{1-1/n}} <(2^{-2/n}-2^{-1/n}\epsilon) \cdot \frac { {\pi}^{1/2} n}{\{ \Gamma(1+n/2)\}^{1/n}},
$$
where $\Gamma_1^{1,1}= \Gamma_1 \cap \overline {\M_1}$.
\label{lem2-1}
\end{lem}

{\bf Proof.}\ Suppose not, then 
$$
\frac{|\Gamma_1^{1,1}|}{|\Omega_1|^{1-1/n}} \ge (2^{-2/n}-2^{-1/n}\epsilon) \cdot \frac { {\pi}^{1/2} n}{\{ \Gamma(1+n/2)\}^{1/n}},
$$
and 
$$
\frac{|\Gamma_1^{1,2}|}{|\Omega_2|^{1-1/n}} \ge (2^{-2/n}-2^{-1/n}\epsilon) \cdot \frac { {\pi}^{1/2} n}{\{ \Gamma(1+n/2)\}^{1/n}},
$$
where $\Gamma_1^{1,2}= \Gamma_1 \cap \overline {\M_2}$. It follows that
$$
|\Gamma_1| \ge  2 (2^{-2/n}-2^{-1/n}\epsilon) \cdot \frac { {\pi}^{1/2} n |\M_1|^{1-1/n}}{\{ \Gamma(1+n/2)\}^{1/n}}=(2^{-1/n}-\epsilon) \cdot \frac { {\pi}^{1/2} n |\M|^{1-1/n}}{\{ \Gamma(1+n/2)\}^{1/n}}.
$$
This contradicts assumption (\ref{A2-1}).

\medskip

\begin{lem} Suppose that $H_\theta$ is a tangent plane of $\Gamma_1 \cap \Gamma_2$ and divides $\M$ into two   domains $\M_1$ and $\M_2$.  If $\overline \M_1$ is the domain which does not contain $\Gamma_1 \cap\Gamma_2$ completely and  $|\M_1| \ge |\M_2|$,  then  on  $\M_2$ the following inequality holds:
$$
\frac{|\Gamma_1^{1,2}|}{|\Omega_2|^{1-1/n}} <(2^{-2/n}-2^{-1/n}\epsilon) \cdot \frac { {\pi}^{1/2} n}{\{ \Gamma(1+n/2)\}^{1/n}},
$$
where $\Gamma_1^{1,2}= \Gamma_1 \cap \overline {\M_2}$.
\label{lem2-2}
\end{lem}

{\bf Proof.}\ Let $|\M_1|=r |\M_2|$ with $r \ge 1$. Using the even reflection of $\M_1$ with respect to $H_\theta$, we have 
$$
\frac {|\Gamma_1^{1,1}|}{|\M_1|^{1-1/n}} \ge 2^{-1/n} \cdot \frac { {\pi}^{1/2} n }{\{ \Gamma(1+n/2)\}^{1/n}}
$$
where $\Gamma_1^{1,1}= \Gamma_1 \cap \overline {\M_1}$. If 
$$
\frac{|\Gamma_1^{1,2}|}{|\Omega_2|^{1-1/n}} \ge (2^{-2/n}-2^{-1/n}\epsilon) \cdot \frac { {\pi}^{1/2} n}{\{ \Gamma(1+n/2)\}^{1/n}},
$$
we have
$$
\begin{array}{rll}
|\Gamma_1|& \ge \big (2^{-1/n}|\M_1|^{1-1/n}+(2^{-2/n}-2^{-1/n}\epsilon)|\M_2|^{1-1/n} \big ) \cdot \frac { {\pi}^{1/2} n}{\{ \Gamma(1+n/2)\}^{1/n}}\\
&=(2^{-1/n}r^{1-1/n}+2^{-2/n}-2^{-1/n}\epsilon)|\M_2|^{1-1/n} \cdot \frac { {\pi}^{1/2} n}{\{ \Gamma(1+n/2)\}^{1/n}}\\
&=2^{-1/n}|\M|^{1-1/n} \cdot \frac {r^{1-1/n}+2^{-1/n}-\epsilon}{(1+r)^{1-1/n}} \cdot \frac { {\pi}^{1/2} n}{\{ \Gamma(1+n/2)\}^{1/n}}\\
&>2^{-1/n}|\M|^{1-1/n}\cdot \frac { {\pi}^{1/2} n}{\{ \Gamma(1+n/2)\}^{1/n}}.
\end{array}
$$
Contradiction to (\ref{A2-1}).

\medskip

We now start the procedure of the reflection symmetrization.

Let $\theta$ be  an angle incommensurable with $\pi$. We use $H_\theta$ to cut through $\M$. There are two different situations.  The  location of $H_\theta$ will be determined differently in these two cases. 

{\bf Case 1.}\ If there is a $H_\theta$ such that $H_\theta$ divides $\M$ into two equal volume domains $\M_1$ and $\M_2$ and $H_\theta\cap \Gamma_2 \ne \O$. By Lemma \ref{lem2-1}, we can assume, without loss of generality, that
$$ 
\frac{|\Gamma_1^{1,1}|}{|\Omega_1|^{1-1/n}} <(2^{-2/n}-2^{-1/n}\epsilon) \cdot \frac { {\pi}^{1/2} n}{\{ \Gamma(1+n/2)\}^{1/n}},
$$
where $\Gamma_1^{1,1}= \Gamma_1 \cap \overline {\M_1}$. Similarly, we define $\Gamma_2^{1,1}= \Gamma_2 \cap \overline {\M_1}$. We then reflect $\M_1$ with respect to $H_\theta$ and obtain $\M_1'$. Let $\M^1=\M_1 \cup \M_1' \cup (H_\theta \cap \M)$, $\Gamma_1^1$= the union of $\Gamma_1^{1,1}$ and its reflection, and  $\Gamma_2^1$= the union of $\Gamma_2^{1,1}$ and its reflection. We have
\bg
\frac{|\Gamma_1^{1}|}{|\Omega^1|^{1-1/n}} <(2^{-1/n}-\epsilon) \cdot \frac { {\pi}^{1/2} n}{\{ \Gamma(1+n/2)\}^{1/n}}.
\label{A2-2}
\ede
Moreover, $\M^1$ has the following two key properties:

{\it Property 1}. If $\{(x_1, ..., x_n) \ : \ x_n=k \} \cap \Gamma_2^1$ has a  strictly concave (with respect to $\M^1$) closed  surface on the hyperplane $\{(x_1, ..., x_n) \ : \ x_n=k \}$,  we denote $M^{n-1}$ as the domain in the hyperplane bounded by the closed surface. Then in $M^{n-1}$ there are  points which are in a neighborhood of a boundary point on $\partial M^{n-1}$,  and are  not in $\M^1$. 

{\it Property 2}. The area of the projection of $\Gamma_2^1$ to  $\{(x_1, ..., x_n) \ : \ x_n=0 \}$  is bounded from below by a fixed constant $C_b>0$  and above by another fixed constant $C_a>0$.

{\it Property 1} follows directly from the assumption that $\Gamma_2$ is concave with respect to $\M$. Since  $|\M^1|=|\M|$ and the height of $\M^1$ is bounded, and $\M$ is a bounded domain, we easily see that {\it Property 2} holds for $\M^1$.

\medskip

{\bf Case 2.}\ If there is no hyperplane which is  parallel to $x_n$-axis and has $\theta$ angle with $x_1$-axis such that it  divides $\M$ into two equal volume domains $\M_1$ and $\M_2$ and has a nonempty intersection with $ \Gamma_2 $. By Lemma \ref{lem2-2}, we know that there is a domain $\M_1$ which is  bounded by $\M$ and a tangent plane $ H_\theta$ of $\Gamma_1 \cap \Gamma_2$, such that $|\M_1| \le |\M|/2$ and
$$
\frac{|\Gamma_1^{1,1}|}{|\Omega_1|^{1-1/n}} <(2^{-2/n}-2^{-1/n}\epsilon) \cdot \frac { {\pi}^{1/2} n}{\{ \Gamma(1+n/2)\}^{1/n}},
$$
where $\Gamma_1^{1,1}= \Gamma_1 \cap \overline {\M_1}$.
 Similarly, we define $\Gamma_2^{1,1}= \Gamma_2 \cap \overline {\M_1}$.  We then reflect $\M_1$ with respect to $H_\theta$ and obtain $\M_1'$. Let  $\M^1=\M_1 \cup \M_1' \cup (H_\theta \cap \M)$, $\Gamma_1^1$= the union of $\Gamma_1^{1,1}$ and its reflection, and  $\Gamma_2^1$= the union of $\Gamma_2^{1,1}$ and its reflection.  It is easy to  check that (\ref{A2-2}) holds for $\M^1$ and $\Gamma_1^1$, and  $\M^1$ also has {\it Property 1} and {\it Property 2}.

\medskip

We next apply $H_{2 \theta}$ on $\M^1$ and obtain $\M^2$ with {\it Property 1} and {\it Property 2}, and so on. Eventually, we obtain a limit domain $\M^{\infty}$ with boundary $\Gamma_1^{\infty}$ and $\Gamma_2^{\infty}$. Due to {\it Property 2}, we know that  $\M^{\infty}$ will not  degenerate into a straight line parallel to $x_n$-axis. From (\ref{A2-2}) we know that  $\M^{\infty}$ will not  degenerate into a $n-1$ domain parallel to the hyperplane $\{(x_1, ..., x_n ) \ : \ x_n=0\}$ neither.   Moreover,
\bg
\frac{|\Gamma_1^{\infty}|}{|\Omega^{\infty}|^{1-1/n}} \le (2^{-1/n}-\epsilon) \cdot \frac { {\pi}^{1/2} n}{\{ \Gamma(1+n/2)\}^{1/n}}.
\label{A2-3}
\ede

On the other hand, it is easy to see that $\Gamma_1^{\infty}$ and $ \Gamma_2^{\infty}$ are revolution surfaces, and $\Gamma_1^{\infty} \cap \Gamma_2^{\infty}$ is a sphere lying on a hyperplane $P_\tau=\{(x_1,...,x_n) \ : \ x_n=\tau\}$ for some constant $\tau$. Due to {\it Property 1}, we know that $\Gamma_2^\infty$ lies on or above the hyperplane $P_\tau$. Let $\M_a^\infty=\{ (x_1,...,x_n) \in \M^\infty \ : \ x_n>\tau\}$, $\M_b^\infty=\{ (x_1,...,x_n) \in \M^\infty \ : \ x_n<\tau\}$, $\Gamma_{1,a}^\infty=\{ (x_1,...,x_n) \in \Gamma_1^\infty \ : \ x_n>\tau\}$, 
and $\Gamma_{1,b}^\infty=\{ (x_1,...,x_n) \in \Gamma_1^\infty \ : \ x_n<\tau\}$.
Due to {\it Property 1}, using the  even reflection with respect to $P_\tau$, we have
$$
{|\Gamma_{1,a}^{\infty}|} \ge 2^{-1/n} \cdot \frac { {\pi}^{1/2} n}{\{ \Gamma(1+n/2)\}^{1/n}}\cdot {|\Omega_a^{\infty}|^{1-1/n}},
$$
and 
$$
{|\Gamma_{1,b}^{\infty}|} \ge 2^{-1/n} \cdot \frac { {\pi}^{1/2} n}{\{ \Gamma(1+n/2)\}^{1/n}} \cdot {|\Omega_b^{\infty}|^{1-1/n}}.
$$
 These two inequalities yields
$$
\frac{|\Gamma_1^{\infty}|}{|\Omega^{\infty}|^{1-1/n}} \ge 2^{-1/n} \cdot \frac { {\pi}^{1/2} n}{\{ \Gamma(1+n/2)\}^{1/n}}.
$$
This contradicts (\ref{A2-3})!  Therefore (\ref{A2-1}) is false for any given $\epsilon<\frac {1+2^{-1/n}-2^{1-1/n}}{2}$. We hereby complete the proof of Theorem \ref{thmA2}.

\section{Equimeasurable rearrangement}

Let $\M \subset \Bbb{R}^n$ be a bounded piecewise smooth domain and $f(x)= f(x_1,..., x_n) \in C^1(\M) \cap C(\overline \M)$ be  a  nonnegative function. 
In this section, we present an integral inequality for $f(x)$ under Schwarz symmetrization. This result was obtained by Duff in \cite{Du}. We present it  here  in a suitable way for our convenience. For completeness, we include all rigorous proofs here.

Let $\mu(t)=meas \{x \in \M \ : \ f(x) >t\}.$  We define a one-dimensional decreasing rearrangement $f^{\#}(t)$ of $f(x)$ as the inverse function of $\mu(t)$, that is: 
\bg
\mu(f^{\#}(t))= t.
\label{A3-1}
\ede
Thus the domain of $f^{\#}(t)$ is $[0, |\M|]$ and the range of $f^{\#}(z)$ is $[\min f(x), \max f(x)].$  

The  basic relation between $d f^{\#}(t)/dt$ and $\nabla f(x)$ is given by the following lemma.

\begin{lem}
For almost every value $f^{\#}(r)$ in $[\min f(x), \max f(x)]$,
$$
\frac 1{ |{f^{\#}}'(r)|} =\sum \int_{f=f^{\#}(r)} \frac {dS}{|\nabla f|},
$$
where ${f^{\#}}'(r)= df^{\#}(t)/dt \mid_{t=r}$, and the integration and summation on the right hand side  run over all components of the level set $f=f^{\#}(r)$ in $\M$ when $\{x \in \M \ : \ f(x)=f^{\#}(r)\}$ consists of a finite number of smooth  surfaces.
\label{lem3-1}
\end{lem}

{\bf Proof.}\ From coarea formula and Sard's Theorem (see, for example, \cite{Ta} or \cite{Evans}), we know that  for almost every $f^{\#}(r)$, $\{x \in \M \ : \ f(x)=f^{\#}(r)\}$ consists of a finite number of smooth surfaces, and
\bg
\frac {d \mu (t)}{dt} \mid_{t=f^{\#}(r)}=-\sum \int_{f=f^{\#}(r)}\frac{ dS}{|\nabla f|}.
\label{A3-2}
\ede

On the other hand, from (\ref{A3-1})  we have
\bg
\frac{d \mu (t)}{dt} \mid_{t=f^{\#}(r)} \cdot \frac {d f^{\#}(r)}{dr} =1.
\label{A3-3}
\ede
Lemma \ref{lem3-1} follows from  (\ref{A3-2}) and (\ref{A3-3}) immediately.

\medskip

\begin{lem} If $\{x \in \M \ : \ f(x)=f^{\#}(r)\}$ consists of a finite number of smooth surfaces, then for any $p > 1$,
$$
(\sum \int_{f=f^{\#}(z)} \frac {dS}{|\nabla f|})^{1-p} \le \frac 1{(S(f))^p}\sum \int_{f=f^{\#}(z)} |\nabla f|^{p-1}dS,
$$
where $S(f)= \sum \int_{f=f^{\#}(z)}dS$.
\label{lem3-2}
\end{lem}

{\bf Proof.} It is equivalent to find the minimum of 
\bg
\sum \int_{f=f^{\#}(z)} g(s)^{1-p}dS
\label{A3-4}
\ede
for positive function $g(s)$  under the constraint: 
$$
\sum \int_{f=f^{\#}(z)} g(s) dS=1.
$$

Consider
$$\sum \int_{f=f^{\#}(z)} \{g(s)^{1-p}+\lambda g(s)\}dS,
$$
where $\lambda$ is a Lagrange multiplier. First variation of $g$ yields the equation
$$
\frac {1-p}{g(s)^p}+\lambda=0,
$$
whence 
$$
g(s)=1/S(f).
$$

Second variation of $g$ is
$$
p(p-1) \sum \int_{f=f^{\#}(z)} \frac {(\delta g(s))^2}{g(s)^{p+1}} dS,
$$
which indicates that (\ref{A3-4}) attains its minimum at $g(s)=1/S(f)$. This yields Lemma \ref{lem3-2}.

\medskip

We are now able  to establish the following key integral inequality.
\begin{prop}
 For any $p > 1$,
\bg
\int_0^{|\M|} |\frac {d f^{\#}(z)}{dz}|^p \cdot (S(f))^pdz \le \int_\M|\nabla f|^p dV,
\label{A3-5}
\ede
where $S(f)= \sum \int_{f=f^{\#}(z)}dS$ for almost every $f^{\#}(z)$.
\label{prop3-1}
\end{prop}

{\bf Proof.}\  From Lemma \ref{lem3-1} and Lemma \ref{lem3-2}, we have
$$
|\frac  {d f^{\#}(z)}{dz}|^{p-1}= (\sum \int_{f=f^{\#}(z)} \frac {dS}{|\nabla f|})^{1-p} \le \frac 1{(S(f))^p}\sum \int_{f=f^{\#}(z)} |\nabla f|^{p-1}dS.
$$
Thus
$$
\begin{array}{rll}
\int_0^{|\M|} |\frac {d f^{\#}(z)}{dz}|^p \cdot (S(f))^pdz &= \int_{z=0}^{z=|\M|}  |\frac {d f^{\#}(z)}{dz}|^{p-1} \cdot (S(f))^p |d f^{\#}(z)|\\
& \le \int_{z=0}^{z=|\M|} \sum \int_{f=f^{\#}(z)}|\nabla f|^{p-1} dS |df^{\#}(z)|\\
&=\int_{\min f(x)}^{\max f(x)} \sum \int_{f=t}|\nabla f|^{p-1} dS dt \\
&=\int_\M|\nabla f|^p dV.
\end{array}
$$
We use the coarea formula in the last equality. This completes the proof of Proposition \ref{prop3-1}.

\medskip

Now, we introduce a radial rearrangement $f^*(x)=f^*(|x|): =f^{\#}(t)$, where
$$
t= \omega_{n-1} |x|^n/n, \ \ \ \ \ \ \ \ \ |x|^2=\sum_{i=1}^n x_i^2,
$$
and $\omega_{n-1}$ is the surface area of the $n-1$-dimensional unit sphere. 
Since $ dt=\omega_{n-1} |x|^{n-1} d|x|$, we have
$$
|\nabla f^*(x)|=| {df^{\#}(t)}/{d t}| \cdot {\omega_{n-1} |x|^{n-1}} =|df^{\#}(t)/dt | \cdot {S(f^*)},
$$
where $S(f^*)= \int_{f^*=f^{\#}(t)} dS$ is the surface area of the level surface $f^*=f^{\#}(t)$, which is a sphere in $\Bbb{R}^n$. 
Combining this with Proposition \ref{prop3-1}, we have
\begin{cor} For any $p>1$,
$$
\int_{\M^*}|\nabla f^*|^p \cdot \frac {(S(f))^p}{(S(f^*))^p} dV \le \int_{\M} |\nabla f|^p dV,
$$
\label{cor3-1}
\end{cor}
where  $S(f)= \int_{f=f^*(x)} dS$ if $\{z \in \M \ : \ f(z)=f^{*}(x)\}$ 
consists of a finite number of smooth surfaces and  $S(f^*)= \int_{f^*=f^*(x)} dS$ .

\section{Proofs of the theorems}

Using the integral inequality in Corollary \ref{cor3-1}, we are able to prove Theorem \ref{thmA3}, \ref{thmA4} and \ref{thmA5} from Theorem \ref{A2}.  Theorem \ref{A1} is simply a corollary of  Theorem \ref{thmA3}.

We need the following  key lemma in our proofs.

\begin{lem}  Let $\M$ be a bounded piecewise smooth domain in $\Bbb{R}^n$ ($n \ge 2$), whose boundary $\partial \M$ consists of one set of disjoint smooth hypersurfaces  $\Gamma_1$   and another  smooth hypersurface $\Gamma_2$, and $\Gamma_2$ is concave with respect to $\M$.  For any  $u \in C^1(\M)\cap C(\overline \M)$ with $u=0$ on any hypersurface in $\Gamma_1$ and $u\ge 0$, we denote $u^*(x)$ as the radial rearrangement of $u(x)$ as being defined in Section 3, and $\M^*$ as the ball with the same volume as $\M$.  Then for any $p>1$,
$$
\int_{\M^*} |\nabla u^*|^p dx \le 2^{p/n} \int_\M |\nabla u|^p dx.
$$
\label{lem4-1}
\end{lem}

{\bf Proof.} \  Due to Corollary \ref{cor3-1} and Sard's Theorem,  we only need to show that for almost every $t$, the level surfaces of  $\{x \in \M \ : \ u=t \}$ and $\{x \in \M^* \ : \ u^*(x)=t\}$ have the following property:
\bg
\frac {S(u)}{S(u^*)} \ge \frac 1{2^{1/n}}.
\label{A4-1}
\ede
Without loss of generality, we can assume that all level surfaces   $\{x \in \M \ : \ u=t \}$ and $\{x \in \M^* \ : \ u^*(x)=t\}$ are regular.

Denote $\M_t= \{x \in \M \ : \ u(x) > t \}$ and $\M^*_t = \{ x \in \M^* \ : \ u^*(x) > t\}.$ We know that $|\M_t|=|\M_t^*|$.

If $t=0$,  by Theorem \ref{thmA2}, we have
$$
\frac {S(u)} {|\M |^{1-1/n}}\ge \frac {|\Gamma_1|} {|\M |^{1-1/n}}  \ge\frac {\pi^{1/2}n}{\{2\Gamma(1+n/2)\}^{1/n}}.
$$
Since for any $t \ge 0$,
\bg
\frac {S(u^*)} {|{\M^*}_t|^{1-1/n}}= \frac {|\partial {\M^*}_t|}{|{\M^*}_t|^{1-1/n}}=\frac {\pi^{1/2}n}{\{\Gamma(1+n/2)\}^{1/n}},
\label{A4-2}
\ede
we have 
$$
\frac {S(u)}{S(u^*)} \ge \frac 1{2^{1/n}},
$$
thus (\ref{A4-1}) holds in this case.

If $t>0$, we write  $\partial \M_t\setminus \Gamma_2=(\cup_{i=1}^k L_i) \cup ( \cup_{j=1}^l \Upsilon_j)$ where $L_i, \ \Upsilon_j$ are  pairwise disjoint smooth surfaces, and satisfy $ \overline L_i \cap \Gamma_2=\O$ and $\overline \Upsilon_j \cap \Gamma_2 \ne \O.$ Since $u=0 $ on $\Gamma_1$, $\partial \M_t \cap \Gamma_1 =\O$. Thus $L_i$ is closed.  Let $  \Phi_i$ be the domain bounded by $L_i$ and $\Psi_j$ be the domain bounded by $\Upsilon_j$ and $\Gamma_2$. 

From the standard isoperimetric inequality, we have
\bg
\frac{|L_i|}{|\Phi_i \cap \M|^{1-1/n}} \ge \frac{|L_i|}{|\Phi_i|^{1-1/n}} \ge  \frac{| \partial B_1|}{|B_1|^{1-1/n}} =\frac {\pi^{1/2}n}{\{\Gamma(1+n/2)\}^{1/n}}, \ \ \ \ \ \ \ \ \mbox{for} \ \ \ i=1,...,k.
\label{A4-3}
\ede
On the other hand, since $\Gamma_2$ is concave we know from Theorem \ref{thmA2}  that
\bg
\frac{|\Upsilon_j|}{|\Psi_j|^{1-1/n}}   \ge \frac {\pi^{1/2}n}{\{2\Gamma(1+n/2)\}^{1/n}}, \ \ \ \ \ \ \ \ \mbox{for} \ \ \ j=1,...,l.
\label{A4-4}
\ede
It follows from (\ref{A4-3}) and (\ref{A4-4}) that
\bg
\begin{array}{rll}
S(u) &= \sum_{i=1}^k |L_i| +\sum_{j=i}^l |\Upsilon_j|\\
&\ge \frac {\pi^{1/2}n}{\{2\Gamma(1+n/2)\}^{1/n}} \cdot( \sum_{i=1}^k |\Phi_i\cap \M |^{1-1/n}+\sum_{j=1}^l |\Psi_j|^{1-1/n}) \\
&  \ge \frac {\pi^{1/2}n}{\{2\Gamma(1+n/2)\}^{1/n}} \cdot( \sum_{i=1}^k |\Phi_i\cap \M|+\sum_{j=1}^l |\Psi_j|)^{1-1/n}\\
& = \frac {\pi^{1/2}n}{\{2\Gamma(1+n/2)\}^{1/n}} |\M_t|^{1-1/n}.
\end{array}
\label{A4-5}
\ede
We find on comparing   (\ref{A4-2}) and (\ref{A4-5}) that (\ref{A4-1}) holds for $t >0$. 
 We hereby complete the proof of Lemma \ref{lem4-1}.

\medskip

\noindent{\bf Proof of Theorem \ref{thmA3}}.

We only need to prove the sharp Sobolev inequality  for $p \in (1, n)$ and for any  $u \in C^1(\M)\cap C(\overline \M)$ with $u=0$ on any hypersurface in $\Gamma_1$ and $u\ge 0$.  
The sharp Sobolev inequality for  $p=1$ follows from taking $p \to 1$.
We denote $u^*(x)$ as the radial rearrangement of $u(x)$ as being defined in Section 3, and $\M^*$ as the ball with the same volume as $\M$.   The equimeasurable property of the rearrangement yields
$$
\int_{\M} |u|^{p_*} = \int_{\M_*} |u^*|^{p_*}.$$
Notice that $u^*=0$ on $\partial \M^*$. It then follows from Lemma \ref{lem4-1} and the sharp Sobolev inequality for $u^* \in W^{1,p}_0(\M_*)$ that 
$$
 \frac {\int_{\M}|\nabla u|^p}{(\int_{\M} |u|^{p_*})^{p/p_*}}\ge  \frac {\int_{\M^*}|\nabla u^*|^p}{2^{p/n}(\int_{\M^*} |u^*|^{p_*})^{p/p_*}}
 \ge \frac 1{2^{p/n } k^p(n,p)}.
$$
This proves inequality (\ref{A8}).

On the other hand, it is well known that
$$
\inf_{v \in C^{1}(\M)\setminus \{0\}, \ v=0 \ on  \ \Gamma_1}  \frac {\int_{\M}|\nabla v|^p}{(\int_{\M} |v|^{p_*})^{p/p_*}} \le  \frac 1{2^{p/n } k^p(n,p)}.
$$ 
This yields the sharpness of the constant $2^{1/n } k(n,p)$.

\medskip

\noindent{\bf Proof of Theorem \ref{thmA4}}.

 For any  $u \in C^1(\M)\cap C(\overline \M)$ with $u=0$ on any hypersurface in $\Gamma_1$ and $u\ge 0$, we denote $u^*(x)$ as the radial rearrangement of $u(x)$ as being defined in Section 3, and $\M^*$ as the ball with the same volume as $\M$.   Notice that $u^*=0$ on $\partial \M^*$. 

 From Lemma \ref{lem4-1} we know that  $||\nabla u^*||_{L^n(\M^*)}^n \le 2$, thus $\int_{\M^*} | 2^{-1/n} \nabla u^*|^n \le 1$. It then follows from the equimeasurable property of the rearrangement and the  Moser-Trudinger inequality for $u^* \in W_0^{1, n}(\M^*)$ that
$$
\int_\M e^{\beta_n  u^{\frac n{n-1}}}= \int_{\M^*} e^{n \omega_{n-1}^{1/(n-1)} (2^{-1/n} u^*)^{n/(n-1)}} \le C_*(|\M|)
$$
for some constant $C_*(|\M|)$ depending only on the volume of $\M$. This completes the proof of Theorem \ref{thmA4}.

\medskip

The assumption of concavity about  $\Gamma_2$ in Theorem \ref{thmA4} is crucial. Without it, Theorem \ref{thmA4} may not be true. Here is

\noindent{\bf A counterexample}.

Let $\Gamma_{i,2}$ be  the curve $y=a_i x^2$ for $-a_i^{-1/3} \le x \le a_i^{-1/3}$ and $\Gamma_{i,1} $ be another curve which connects $(-a_i^{-1/3}, a_i^{1/3})$ with $(a_i^{-1/3}, a_i^{1/3})$ and is above the straight line $y=a_i^{1/3}$ such that the area of the region $\M_i$ bounded by $\Gamma_{i,1}$ and $\Gamma_{i,2}$ is $1$, where $a_i \to \infty$ as $i \to \infty$.

For $\E <e^{-a_i^2}$, we define
$$
u_i=
\left\{
\begin{array}{rll}
&\frac {\sqrt 2}{\sqrt{2 \pi}} \cdot (\ln \frac 1 \E)^{1/2}, \ \ \ \ & 0 \le |x|\le \E,\\
&\frac {\sqrt 2}{\sqrt{2 \pi}} \cdot \frac {\ln \frac 1{|x|}}{(\ln \frac 1 \E)^{1/2}}, \ \ \ \ & \E\le |x| \le 1,\\
&0, \ \ & \mbox{otherwise}.
\end{array}
\right.
$$
Claim:
$$
\int_{\M_i} e^{2 \pi u_i^2/||\nabla u_i||^2_{L^2(\M_i)} } \to \infty \ \ \ \ \ \ \ \mbox{as} \ \ \ i \to \infty.
$$

We first convert $(x,y)$ into polar coordinates $(r, \theta)$. It is easy to see that on $\Gamma_{i,2}$  there is a $ 0< \tau_0<1/100$ such that
for $\tau_0 a_i^{-1} \le  r \le 1$, $\theta>\tau_0/2$ in the first quadrant, and  $\pi-\theta>\tau_0/2$ in the second quadrant.

Thus
$$
\int_{\M_i} |\nabla u_i|^2 dz  \le \int_\E^{\tau_0 a_i^{-1}} \frac {\pi}{\pi r \ln \frac 1 \E} dr + \int_{\tau_0 a_i^{-1}}^1 \frac {\pi-\tau_0}{\pi r \ln \frac 1 \E} dr  = 1- \frac {{\tau_0} \ln \tau_0^{-1}a_i}{\pi \ln \frac 1 \E}.
$$

For $|x|< \E$,
$$
\frac {2 \pi u_i^2}{||\nabla u_i||_2^2} =\frac {2 \ln 1/ \E}{1-\frac {{\tau_0} \ln \tau_0^{-1}a_i}{\pi \ln \frac 1 \E} } \ge  2  (1+\frac {{\tau_0} \ln \tau_0^{-1}a_i}{\pi \ln \frac 1 \E} ) \cdot \ln \frac 1 \E.
$$
Therefore
$$
\begin{array}{rll}
\int_{\M_i} e^{2 \pi u_i^2/||\nabla u_i||^2_{L^2(\M_i)} }& \ge \int_{\{(x, y) \in \M_i \ : \ x^2+y^2 <\E^2\}} e^{ 2 \pi u_i^2/||u_i||_2^2} \ge \frac {\pi \E^2 (1+o_i(1))}{\E^{2+\frac {2 \tau_0 \ln a_i \tau_0^{-1}}{\pi \ln 1/\E}}}  \\
&  =\pi  (1+o_i(1)) e^{\frac{2 \tau_0 \ln a_i \tau_0^{-1}}{\pi}} \to \infty \ \ \ \ \ \ \ \mbox{as} \ \ \ i \to \infty,
\end{array}
$$
where  $o_i(1) \to 0$ as $i \to \infty$ and  we use $\E^{-1/(\ln1/\E)} =e$.

\medskip

\noindent{\bf Proof of Theorem \ref{thmA5}}.

 Assume that $u \in H^1(\M)\setminus \{0\}$ with $u=0$ on $\Gamma_1$.  Let $u^*(x)$ be the radial rearrangement of $|u(x)|$ as being defined in Section 3, and $\M^*$ be the ball with the same volume as $\M$. Notice that $u^*=0$ on $\partial \M^*$. Using Lemma \ref{lem4-1}, we have:
$$
\int_{\M^*} |\nabla u^*|^2 dx \le 2^{2/n} \int_{\M}|\nabla |u||^2 dx.
$$
This yields
\bg
\frac {\int_\M |\nabla u|^2 dx}{\int_\M u^2 dx}= \frac {\int_\M |\nabla |u||^2 dx}{\int_\M u^2 dx}\ge \frac {\int_{\M^*} |\nabla u^*|^2}{2^{2/n} \int_{\M^*} (u^*)^2 dx} \ge \frac{\lambda_1(\M_*)}{2^{2/n}},
\label{A4-6}
\ede
where $\lambda_1(\M^*)$ is the first eigenvalue of $-\Delta$ with Dirichlet boundary condition in $\M^*$. On the other hand, using an even reflection and a  rescaling, we  find that
\bg
\frac{\lambda_1(\M^*)}{2^{2/n}} = \inf_{u\in H^1(B_r^+) \setminus \{0\}, u=0 \ \mbox{on}\ \partial B_r^+\setminus \{y=0\}} \frac{\int_{B_r^+} |\nabla u|^2 dx}{\int_{B_r^+} u^2 dx},
\label{A4-7}
\ede
where $B_r^+$ is the upper half ball with center at $\{0\}$ and  the same volume as $\M^*$. From (\ref{A4-6}) and (\ref{A4-7}) we obtain  Theorem \ref{thmA5}.

\end{document}